**Oleksandra Kukharenko,**
**Taras Shevchenko National University of Kyiv, Ukraine**
*akukharenko@ukr.net*

# SOLUTION OF THE CAUCHY PROBLEM FOR OSCILLATORY EQUATION WITH DELAY

*В даній роботі розглянуто диференціальне рівняння другого порядку з одним сталим запізнюванням. Використовуючи метод кроків, отримано аналітичний розв'язок цього рівняння.*

*Ordinary differential equations of the second order with one constant delay are considered in this paper. An analytical representation of the solution is obtained using the method of steps.*

**Introduction**

Many mathematical models of oscillating processes are described by linear differential equations of the second order [1,2]. If a prehistory has to be considered, what is typical for models in biology, economy, dynamics of populations, delay differential equations [3-6] are used. In the paper [7] the equations and systems of differential equations with pure delay have been considered. Special functions have been introduced and a solution of the Cauchy problem has been received.

It is known, a solution of linear inhomogeneous differential equation of the second order

$$x''(t) + \omega^2 x(t) = f(t), \ t > 0, \tag{1.1}$$

with initial condition $x(0) = x_0$, can be written as

$$x(t) = x_0 \cos \omega t + \frac{x'_0}{\omega} \sin \omega t + \int_0^t \sin \omega(t-s) f(s) ds. \tag{1.2}$$

For the delay equation

$$x''(t) + \omega^2 x(t-\tau) = f(t), \ t > 0, \ \tau > 0,$$

with initial conditions $x(t) = \varphi(t)$, $-\tau \leq t \leq 0$, an analogical representation has the next form

$$x(t) = \varphi(-\tau)\cos_\tau(\omega,t) + \frac{1}{\omega}\dot{\varphi}(-\tau)\sin_\tau(\omega,t) + \frac{1}{\omega}\int_{-\tau}^{0} \sin_\tau(\omega,t-\tau-\xi)\ddot{\varphi}(\xi)d\xi + \int_0^t \sin_\tau(\omega,t-\tau-\xi)f(\xi)d\xi,$$

where $\cos_\tau(\omega,t)$, $\sin_\tau(\omega,t)$ are the special functions named the delay cosine and sine functions [7]. They can be presented as follows.

The delay cosine function $\cos_\tau(\omega,t)$ is a function which can be written as

$$\cos_\tau(\omega,t) = \begin{cases} 0, & -\infty < t < -\tau, \\ 1, & -\tau \leq t < 0, \\ 1 - \omega^2 \dfrac{t^2}{2!}, & 0 \leq t < \tau, \\ \ldots & \ldots \\ 1 - \omega^2 \dfrac{t^2}{2!} + \omega^4 \dfrac{(t-\tau)^4}{4!} + \cdots + (-1)^k \omega^{2k} \dfrac{[t-(k-1)\tau]^{2k}}{(2k)!}, & (k-1)\tau \leq t < k\tau \end{cases}$$

a $2k$-degree polynomial on intervals $(k-1)\tau \leq t < k\tau$ merged in points $t = k\tau, \ k = 0,1,2,\ldots$.

The delay sine function $\sin_\tau(\omega,t)$ is a function which can be written as

$$\sin_\tau(\omega,t) = \begin{cases} 0, & -\infty < t < -\tau, \\ \omega(t+\tau), & -\tau \leq t < 0, \\ \omega(t+\tau) - \omega^3 \dfrac{t^3}{3!}, & 0 \leq t < \tau, \\ \ldots & \ldots \\ \omega(t+\tau) - \omega^3 \dfrac{t^3}{3!} + \cdots + (-1)^k \omega^{2k+1} \dfrac{[t-(k-1)\tau]^{2k+1}}{(2k+1)!}, & (k-1)\tau \leq t < k\tau \end{cases}$$

a $(2k+1)$-degree polynomial on intervals $(k-1)\tau \leq t < k\tau$ merged in points $t = k\tau, \ k = 0,1,2,\ldots$.

**1. Homogeneous Equation**

In the presented paper a linear inhomogeneous equation with one constant delay is considered

$$x''(t) + \omega_1^2 x(t) + \omega_2^2 x(t-\tau) = f(t), \ t \geq 0, \ \tau > 0, \tag{1.3}$$

with an initial condition

$$x(t) \equiv \varphi(t), \quad -\tau \le t \le 0. \tag{1.4}$$

Here $\varphi(t)$ is an arbitrary, continuously differentiable function on an initial interval $-\tau \le t \le 0$. By a solution of the Cauchy problem (1.4) for the equation (1.3) we mean a function $x(t)$, which is defined for $t \ge -\tau$, continuously differentiable for $t \ge 0$ and identically satisfies (1.3) with initial conditions (1.4). It is known, the solution of the Cauchy problem (1.4) for the inhomogeneous equation (1.3) can be presented in the form of sum

$$x(t) = x^0(t) + \overline{x}(t),$$

where $x^0(t)$ is a solution of homogeneous equation

$$x''(t) + \omega_1^2 x(t) + \omega_2^2 x(t-\tau) = 0, \quad t > 0, \quad \tau > 0, \tag{1.5}$$

with initial conditions (1.4). And $\overline{x}(t)$ is a solution of the inhomogeneous equation (1.3) with zero initial conditions $\overline{x}(t) \equiv 0$, $\overline{x}'(t) \equiv 0$, $-\tau \le t \le 0$.

Let $x_1(t)$ is a solution of the homogeneous differential equation with initial conditions

$$x_1(t) \equiv 1, \quad x_1'(t) \equiv 0, \quad -\tau \le t \le 0, \quad x_1(t) \equiv 0, \quad x_1'(t) \equiv 0, \quad t < -\tau, \tag{1.6}$$

and $x_2(t)$ is a solution of the homogeneous differential equation (1.5) with initial conditions

$$x_2(t) \equiv t+\tau, \quad x_2'(t) \equiv 1, \quad -\tau \le t \le 0, \quad x_1(t) \equiv 0, \quad x_1'(t) \equiv 0, \quad t < -\tau. \tag{1.7}$$

Preliminary we will formulate the statement about representation of solution of the Cauchy problem (1.4) for the homogeneous equation (1.5), using functions $x_1(t)$, $x_2(t)$ with initial conditions (1.6), (1.7). The obtained representation is similar to the one formulated at [4].

**Theorem 1.1.** A solution of the homogeneous equation (1.5) with initial conditions (1.4) can be written as

$$x(t) = \varphi(-\tau)x_1(t) + \varphi'(-\tau)x_2(t) + \int_{-\tau}^{0} x_2(t-\tau-s)\varphi''(s)ds, \tag{1.8}$$

where $x_1(t)$ is the solution of the equation (1.5) with the conditions (1.6), $x_2(t)$ is the solution of (1.5) with the conditions (1.7).

*Proof.* We seek solutions of (1.5) of the form

$$x(t) = ax_1(t) + bx_2(t) + \int_{-\tau}^{0} x_2(t-\tau-s)c''(s)ds, \tag{1.9}$$

where $a$, $b$ are arbitrary constants and $c(t)$ is a twice continuously differentiable function. Substituting (1.9) into the equation (1.5), we obtain

$$\frac{d^2}{dt^2}\left[ax_1(t) + bx_2(t) + \int_{-\tau}^{0} x_2(t-\tau-s)c(s)ds\right] + \omega_1^2\left[ax_1(t) + bx_2(t) + \int_{-\tau}^{0} x_2(t-\tau-s)c''(s)ds\right] +$$

$$+ \omega_2^2\left[ax_1(t-\tau) + bx_2(t-\tau) + \int_{-\tau}^{0} x_2(t-2\tau-\xi)c''(s)ds\right] = 0.$$

We rewrite the obtained expression and get

$$a\left[\ddot{x}_1(t) + \omega_1^2 x_1(t) + \omega_2^2 x_1(t-\tau)\right] + b\left[\ddot{x}_2(t) + \omega_1^2 x_2(t) + \omega_2^2 x_2(t-\tau)\right] +$$

$$+ \int_{-\tau}^{0}\left[\ddot{x}_2(t-\tau-s) + \omega_1^2 x_2(t-\tau-s) + \omega_2^2 x_2(t-2\tau-s)\right]c''(s)ds \equiv 0.$$

By the definition, $x_1(t)$ and $x_2(t)$ are solutions of the homogeneous equation, so expressions in square brackets are identically equal to zero, and the obtained identity shows that (1.9) is a solution of the equation (1.5) for the arbitrary constants $a$, $b$ and the function $c(t)$.

We will show that for $a = \varphi(-\tau)$, $b = \varphi'(-\tau)$, $c(t) = \varphi(t)$ initial conditions are satisfied

$$x(t) = \varphi(t), \quad x'(t) = \varphi'(t), \quad -\tau \le t \le 0,$$

i.e. the solution of the Cauchy problem (1.4) for the homogeneous equation (1.5) has the form (1.8).

Substituting initial conditions (1.4) into (1.8), we obtain

$$\varphi(t) = \varphi(-\tau)x_1(t) + \varphi'(-\tau)x_2(t) + \int_{-\tau}^{0} x_2(t-\tau-s)\varphi''(s)ds, \quad -\tau \le t \le 0.$$

We make the change of variables $t - \tau - s = \xi$ in the integral, so that

$$\varphi(t) = \varphi(-\tau)x_1(t) + \varphi'(-\tau)x_2(t) + \int_{t-\tau}^{t} x_2(\xi)\varphi''(t-\tau-\xi)d\xi, \quad -\tau \leq t \leq 0.$$

The integral is divided into two

$$\varphi(t) = \varphi(-\tau)x_1(t) + \varphi'(-\tau)x_2(t) + \int_{t-\tau}^{-\tau} x_2(\xi)\varphi''(t-\tau-\xi)d\xi + \int_{-\tau}^{t} x_2(\xi)\varphi''(t-\tau-\xi)d\xi.$$

Considering that for $t-\tau \leq \xi \leq -\tau$ $x_1(\xi) \equiv 0$, the first integral vanishes, and

$$\varphi(t) = \varphi(-\tau)x_1(t) + \varphi'(-\tau)x_2(t) + \int_{-\tau}^{t} x_2(\xi)\varphi''(t-\tau-\xi)d\xi, \quad -\tau \leq t \leq 0.$$

The integral is taken by parts

$$\varphi(t) = \varphi(-\tau)x_1(t) + \varphi'(-\tau)x_2(t) - x_2(\xi)\varphi'(t-\tau-\xi)\Big|_{\xi=-\tau}^{\xi=t} + \int_{-\tau}^{t} x_2'(\xi)\varphi'(t-\tau-\xi)d\xi =$$

$$= \varphi(-\tau)x_1(t) + x_2(-\tau)\varphi'(t) + \int_{-\tau}^{t} x_2'(\xi)\varphi'(t-\tau-\xi)d\xi.$$

$x_1(t) \equiv 1$, $x_2(-\tau) = 0$, $x_2'(t) \equiv 1$ in the interval $-\tau \leq t \leq 0$ and we have

$$\varphi(t) = \varphi(-\tau) + \int_{-\tau}^{t} \varphi'(t-\tau-\xi)d\xi = \varphi(-\tau) - \varphi(t-\tau-\xi)\Big|_{\xi=-\tau}^{\xi=t} = \varphi(t).$$

Thus the first initial condition is fulfilled. We will show that for the expression (1.8) the second initial condition is fulfilled too. Differentiating (1.8), we obtain

$$x'(t) = \varphi(-\tau)x_1'(t) + \varphi'(-\tau)x_2'(t) + \int_{-\tau}^{0} x_2'(t-\tau-s)\varphi''(s)ds.$$

After the change of variables in the integral, we have

$$x'(t) = \varphi(-\tau)x_1'(t) + \varphi'(-\tau)x_2'(t) + \int_{t-\tau}^{t} x_2'(\xi)\varphi''(t-\tau-\xi)d\xi.$$

Considering that on the initial interval $x_1'(t) \equiv 0$, $x_2'(t) \equiv 1$, we get

$$x'(t) = \varphi'(-\tau) + \int_{-\tau}^{t} \varphi''(t-\tau-\xi)d\xi = \varphi'(-\tau) - \varphi'(t-\tau-\xi)\Big|_{\xi=-\tau}^{\xi=t} = \varphi'(t).$$

The second initial condition is fulfilled.

Thus, to obtain a solution of the Cauchy problem it is necessary to have functions $x_1(t)$, $x_2(t)$, which are solutions of homogeneous equation with special (unit) initial conditions.

**Theorem 1.2.** Solutions $x_1(t)$, $x_2(t)$ of the equations (1.5) with initial conditions (1.6), (1.7) can be presented in the following form

$$x_1(t) = \begin{cases} x_{10}(t), & -\tau \leq t < 0, \\ x_{11}(t), & 0 \leq t < \tau, \\ \quad \ldots \\ x_{1k}(t), & (k-1)\tau \leq t < k\tau, \end{cases} \quad x_2(t) = \begin{cases} x_{20}(t), & -\tau \leq t < 0, \\ x_{21}(t), & 0 \leq t < \tau, \\ \quad \ldots \\ x_{2k}(t), & (k-1)\tau \leq t < k\tau, \end{cases} \quad (1.10)$$

where

$$x_{1k}(t) = x_{1k-1}((k-1)\tau)\cos\omega_1(t-(k-1)\tau) + \frac{x_{1k-1}'((k-1)\tau)}{\omega_1}\sin\omega_1(t-(k-1)\tau) - \omega_2^2 \int_{(k-1)\tau}^{t} \sin\omega_1(t-s)x_{1k-1}(s)ds,$$

$$(k-1)\tau \leq t < k\tau, \quad x_{10}(t) = 1, \quad x_{10}'(t) = 0,$$

$$x_{2k}(t) = x_{2k-1}((k-1)\tau)\cos\omega_1(t-(k-1)\tau) + \frac{x_{2k-1}'((k-1)\tau)}{\omega_1}\sin\omega_1(t-(k-1)\tau) - \omega_2^2 \int_{(k-1)\tau}^{t} \sin\omega_1(t-s)x_{2k-1}(s)ds,$$

$$(k-1)\tau \leq t < k\tau, \quad x_{20}(t) = t+\tau, \quad x_{20}'(t) = 1.$$

(1.11)

*Proof.* To obtain solutions $x_1(t)$, $x_2(t)$ of the equation (1.5) we will use a method of steps [4]. As it was shown above, the solution of the linear inhomogeneous equation without delay (1.1) is (1.2).

Firstly, we consider the first interval $T_1 : 0 \leq t < \tau$. In this interval (1.5) is a linear inhomogeneous differential equation without delay of the following form

$$x''(t) + \omega_1^2 x(t) = -\omega_2^2 x_{10}(t-\tau) \tag{1.12}$$

with the right-hand side presented as a known function $x_{10}(t) = 1$, $x'_{10}(t) = 0$. As it was shown in (1.2), a solution of (1.12) on the interval $T_1 : 0 \leq t < \tau$ will be

$$x_{11}(t) = x_{10}(0)\cos\omega_1 t + \frac{x'_{10}(0)}{\omega_1}\sin\omega_1 t - \omega_2^2 \int_0^t \sin\omega_1(t-s) x_{10}(s)\,ds. \tag{1.13}$$

Next, we consider the second interval $T_2 : \tau \leq t < 2\tau$. In this interval the equation (1.5) has the next form

$$x''(t) + \omega_1^2 x(t) = -\omega_2^2 x_{11}(t-\tau). \tag{1.14}$$

a solution of (1.14) in the interval $T_2 : \tau \leq t < 2\tau$ has a similar representation

$$x_{12}(t) = x_{11}(\tau)\cos\omega_1(t-\tau) + \frac{x'_{11}(\tau)}{\omega_1}\sin\omega_1(t-\tau) - \omega_2^2 \int_\tau^t \sin\omega_1(t-s) x_{11}(s)\,ds. \tag{1.13}$$

Continuing process, we obtain solutions of the homogeneous equation (1.5) on intervals $T_k : (k-1)\tau \leq t < k\tau$:

$$x_{1k}(t) = x_{1k-1}((k-1)\tau)\cos\omega_1(t-(k-1)\tau) + \frac{x'_{1k-1}((k-1)\tau)}{\omega_1}\sin\omega_1(t-(k-1)\tau) - \omega_2^2 \int_{(k-1)\tau}^t \sin\omega_1(t-s) x_{1k-1}(s)\,ds.$$

Thus, $x_1(t)$ has the following form

$$x_1(t) = \begin{cases} x_{10}(t), & -\tau \leq t < 0, \\ x_{11}(t), & 0 \leq t < \tau, \\ \ldots & . \\ x_{1k}(t) & (k-1)\tau \leq t < k\tau. \end{cases}$$

And the solution $x_2(t)$ of (1.5) with initial conditions $x_{20}(t) \equiv t + \tau$, $x'_{20}(t) \equiv 1$, $-\tau \leq t \leq 0$ can be presented in the similar form

$$x_2(t) = \begin{cases} x_{20}(t), & -\tau \leq t < 0, \\ x_{21}(t), & 0 \leq t < \tau, \\ \ldots & . \\ x_{2k}(t) & (k-1)\tau \leq t < k\tau, \end{cases}$$

$$x_{2k}(t) = x_{2k-1}((k-1)\tau)\cos\omega_1(t-(k-1)\tau) + \frac{x'_{2k-1}((k-1)\tau)}{\omega_1}\sin\omega_1(t-(k-1)\tau) - \omega_2^2 \int_{(k-1)\tau}^t \sin\omega_1(t-s) x_{2k-1}(s)\,ds.$$

## 2. Inhomogeneous Equation

Now we will consider the linear inhomogeneous equation with one constant delay (1.3)

$$x''(t) + \omega_1^2 x(t) + \omega_2^2 x(t-\tau) = f(t), \quad t \geq 0, \quad \tau > 0,$$

with zero initial conditions

$$x(t) \equiv 0, \quad x'(t) = 0, \quad -\tau \leq t \leq 0. \tag{2.1}$$

**Theorem 2.1.** A solution of the delay inhomogeneous equation (1.3) with zero initial conditions (2.1) is

$$\overline{x}(t) = \int_0^t x_1(t-s-\tau) f(s)\,ds,$$

where $x_1(t)$ is a solution of the homogeneous equation with initial conditions $x_1(t) \equiv 1$, $x'_1(t) \equiv 0$, $-\tau \leq t \leq 0$.

*Proof.* To proof statements of the theorem we will use the method of variation of parameters. We seek a particular solution $\overline{x}(t)$ of (2.1) in the next form

$$\overline{x}(t) = \int_0^t x_1(t-\tau-s) c(s)\,ds,$$

where $c(t)$ is an unknown function. Substituting this expression into (1.3), we obtain

$$\frac{d}{dt}\left\{\int_0^t x_1(t-\tau-s) c(s)\,ds\right\} + \omega_1^2 \left\{\int_0^t x_1(t-\tau-s) c(s)\,ds\right\} + \omega_2^2 \left\{\int_0^{t-\tau} x_1(t-2\tau-s) c(s)\,ds\right\} = f(t).$$

After calculating derivatives, we have

$$x_1(-\tau)c(t) + \int_0^t x_1'(t-\tau-s)c(s)ds + \omega_1^2\left\{\int_0^t x_1(t-\tau-s)c(s)ds\right\} + \omega_2^2\left\{\int_0^{t-\tau} x_1(t-2\tau-s)c(s)ds\right\} = f(t).$$

We rewrite the obtained expression in the following form

$$x_1(-\tau)c(t) + \int_0^t \{x_1'(t-\tau-s) + \omega_1^2 x_1(t-\tau-s) + \omega_1^2 x_2(t-2\tau-s)\}c(s)ds - \omega_2^2\left\{\int_{t-\tau}^t x_1(t-2\tau-s)c(s)ds\right\} = f(t),$$

and make the change of variables $t - 2\tau - s = \xi$ in the second integral

$$\int_{t-\tau}^t x_1(t-2\tau-s)c(s)ds = \int_{-\tau}^{-2\tau} x_1(\xi)c(t-2\tau-\xi) = 0.$$

Considering that $x_1(-\tau) = 1$ and $x_1(t)$ is a solution of the homogeneous equation, we have

$$c(t) = f(t),$$

i.e. the statements of the theorem 2.1.

Using obtained result we can formulate the following statement.

**Theorem 2.2.** A solution of the delay inhomogeneous equation (1.3) with zero initial conditions (1.4) is

$$x(t) = \varphi(-\tau)x_1(t) + \varphi'(-\tau)x_2(t) + \int_{-\tau}^0 x_2(t-\tau-s)\varphi''(s)ds + \int_0^t x_1(t-s-\tau)f(s)ds,$$

where $x_1(t)$, $x_2(t)$ are defined in (1.10), (1.11).